\documentclass[10pt,a4paper]{article}
\usepackage{amsmath,amssymb}
\usepackage{amsthm}
\usepackage{graphicx}
\usepackage[all]{xy}

\newtheorem{thm}{Theorem}[section]
\newtheorem{lemma}{Lemma}[section]

\newtheorem{remark}{Remark}[section]
\newtheorem{ex}{Example}[section]

\newcommand{\N}{\mathbb N}
\newcommand{\Z}{\mathbb Z}
\newcommand{\R}{\mathbb R}
\newcommand{\Sp}{\mathbb S}

\newcommand{\C}{\mathbb C}

\newcommand{\xJ}{\mathfrak{J}}
\newcommand{\xj}{\mathfrak{j}}
\newcommand{\xp}{\mathfrak{p}}
\newcommand{\lagl}{\mathfrak{gl}}
\newcommand{\lag}{\mathfrak{g}}
\newcommand{\lap}{\mathfrak{p}}
\newcommand{\laso}{\mathfrak{so}}

\newcommand{\gr}{\mathfrak{gr}}
\newcommand{\ms}{\mathfrak{M}}

\newcommand{\dg}{\textbf{deg}}

\newcommand{\SO}{\mathrm{SO}}
\newcommand{\G}{\mathrm G}
\newcommand{\gP}{\mathrm P}
\newcommand{\GL}{\mathrm{GL}}

\newcommand{\Spin}{\mathrm{Spin}}

\newcommand{\wdeg}{\mathfrak{wd}}
\newcommand{\fs}{\mathrm{R}}

\newcommand{\mU}{\mathbb U}
\newcommand{\mV}{\mathbb V}
\newcommand{\mE}{\mathbb E}
\newcommand{\mF}{\mathbb F}
\newcommand{\mW}{\mathbb W}
\newcommand{\tab}{\mathbb A}

\newcommand{\cA}{\mathcal A}
\newcommand{\aff}{\mathcal A}
\newcommand{\cC}{\mathcal C}

\newcommand{\cV}{\mathcal V}

\newcommand{\cU}{\mathcal U}
\newcommand{\cG}{\mathcal G}
\newcommand{\cJ}{\mathcal J}
\newcommand{\cL}{\mathcal L}

\newcommand{\ra}{\rightarrow}

\title{$k$-Dirac operator and the Cartan-K\"ahler theorem for weighted differential operators}
\author{Tom\'a\v s Sala\v c}

\begin{document}
\maketitle
\author

\begin{abstract}
The $k$-Dirac operator is a first order differential operator which is natural to a particular class of parabolic geometries which include the Lie contact structures. A natural task is to understand the set of local null solutions of the operator at a given point. We will show that this set has a very nice and simple structure, namely we will show that there is a submanifold passing through the point such that any section defined on the submanifold extends locally to an unique null solution of the operator. This result also indicates that these parabolic geometries are naturally associated to a certain constant coefficient operator which has been studied in Clifford analysis and this is the original motivation for this paper. In order to prove  the claim about the set of initial conditions for the $k$-Dirac operator we will adapt some parts of the  theory of exterior differential systems and the Cartan-K\"ahler theorem to the setting of differential operators which are natural to geometric structures that are equipped with a filtration of the tangent bundle.\\


\end{abstract}

\section{Introduction}

We will first recall the constant coefficient differential operator from the abstract that has been intensively studied in Clifford analysis (see \cite{CSSS} or \cite{SSSL}). It will be then more natural to explain the motivation for this paper. This operator  generalizes  the $k$-Cauchy-Riemann operator just as the Dirac operator can be viewed as  a generalization of the Cauchy-Riemann operator. 
 
Let $\{\varepsilon_1,\ldots,\varepsilon_{n+1}\}$  be the standard basis  of $\R^{n+1}$, $g$ the standard inner product, $\Sp$ be the complex space of spinors (see Section 6 in \cite{GW}) of the complexified Clifford algebra of $(\R^{n+1},g)$ and $M(n+1,k,\R)$ be the space of matrices of size $(n+1)\times k$. We assume throughout this paper that $k\ge2$ and $n\ge2$. Let $\cC^\infty(M(n+1,k,\R),\mV)$ be the space of all smooth functions on $M(n+1,k,\R)$ with values in some vector space $\mV$. We define a differential operator 
\begin{eqnarray}\label{k-Dirac operator}
E:\cC^\infty(M(n+1,k,\R),\Sp)\ra\cC^\infty(M(n+1,k,\R),\R^k\otimes\Sp)\\
E\psi:=(E_1\psi,\dots,E_k\psi)\ \mathrm{where} \ E_i\psi:=\sum_{\alpha=1}^{n+1}\varepsilon_\alpha.\partial_{x_{\alpha i}}\psi\nonumber
\end{eqnarray}
and $x_{\alpha i}$ are the usual matrix coefficients on $M(n+1,k,\R)$ and the dot denotes the Clifford multiplication. We will call this operator a $k$-\textit{Dirac operator in the flat setting}\footnote{Note that this operator is  called in  \cite{CSSS} and \cite{SSSL} the $k$-Dirac operator but we choose this longer name in order to distinguish it from the $k$-Dirac operator which lives in the world of  parabolic geometries.} and a solution of $E\psi=0$ a \textit{monogenic function}. It is not hard to see that any monogenic function $\psi$ is uniquely determined by its restriction $\psi|_{M(n,k,\R)}$ to the affine subset $M(n,k,\R)=\{x_{n+1,1}=\ldots=x_{n+1,k}=0\}$. Moreover, the restriction has to satisfy the following equations:
\begin{equation}\label{second order equation}
[ \hat E_i, \hat E_j](\psi|_{M(n,k,\R)})=0 ;  \ i,j=1,\dots,k
\end{equation}
where $\hat E_i:=\sum_{\alpha=1}^n\varepsilon_\alpha.\partial_{x_{\alpha i}}$. The equations (\ref{second order equation}) can be easily derived from the fact that the coordinate vector fields commute. It follows that not every $\Sp$-valued function on  $M(n,k,\R)$ is the restriction of a (necessarily unique) monogenic function and it is easy to see  that this implies  (actually it is equivalent to  the fact)  that $E$ is \textit{not} an involutive operator (see   \cite{S} that the first prolongation of the $2$-Dirac operator is involutive and it is probably not known which prolongation of the $k$-Dirac operator with $k>2$ is involutive) in the sense of the Cartan-K\"ahler theorem (see \cite{Br}). 

Recall that the involutivity of the system of PDE's can be most easily checked by the well known Cartan's test (see  \cite{Br}). Let us denote by $M_i$ the space of all monogenic functions whose components are homogeneous polynomials of degree $i$.  Now as $E$ is a first order constant coeffiecient differential operator without zero order part, it follows that a real analytic function  is monogenic\footnote{Actually any monogenic function has to be real analytic as components are harmonic functions.} iff its $i$-th homogeneous component belongs to $M_i$. Moreover, as any coordinate vector field is an infinitesimal symmetry of $E$, differentiation by such a vector field $X$ induces a linear map $M_i\ra M_{i-1}$ which we also denote by $X$.  A choice of an ordered basis $\{X_1,\dots,X_{(n+1)k}\}$ of the vector space of coordinate vector fields  induces filtration  $\{M_i^j:j=1,\dots,(n+1)k\}$ of $M_i$  where $M_i^j:=\{\psi\in M_i:X_1\psi=\dots X_j\psi=0\}$. Notice that since coordinate vector fields commute, the map $X:M_i\ra M_{i-1}$ is compatible with these filtrations. Now it is easy to show that $\dim M_{i+1}\le \dim M_i+ \sum_{j=1}^{k(n+1)}\dim M_i^j$. If the equality holds then we know by the Cartan's test  that the system is involutive and we get an explicit characterization of the set of initial conditions. 

There is a generalization of the Cartan's test for differential operators with real analytic coefficients (see Chapter 9.1 in \cite{Br}). Here one looks at the symbol of the operator at a fixed point and defines a filtration of the kernel of the symbol map. The filtration comes from a choice of an ordered basis of the cotangent space at the point. From this point of view, one simply regards the operator at the given point as a homogeneous\footnote{Meaning homogeneous in the degree of derivatives.} constant coefficient differential operator and uses the Cartan's test we explained above.

In the paper we will consider the differential operator which is called in \cite{S} and \cite{Sa}  the $k$-\textit{Dirac operator in the parabolic setting}. We will denote this operator by $D$ and we will call it for simplicity the $k$-\textit{Dirac operator} and a  solution of $D\Psi=0$ a \textit{monogenic spinor}. The operator $D$ is a  differential operator which is natural to a particular class of parabolic geometries (see Section \ref{section geometric structure}). The geometry is $|2|$-graded which means that any manifold with this geometric structure has  a  bracket generating distribution such that the Levi form has a special algebraic type. If $k=2$, then this is a contact structure  known as the Lie contact structure, see Section 4.2.5 in \cite{CS}. As $D$ is a first order operator, it can be shown (see \cite{SS}) that $D$ is given by differentiating only in the directions which belong to the canonical distribution. This means that $D$ is not just a differential operator of  order one but more 
is true, it is a differential operator of the weighted order one (see Section \ref{section weighted jets} for the definition). 

We will consider $D$ on the homogeneous space $M$ of the parabolic geometry. Given any two points $x,x'\in M$ there is a symmetry $\Phi$ of the geometric structure  which maps $x$ to $x'$. By definition, any natural operator commutes with the induced pullback map $\Phi^\ast$ and so $\Phi^\ast$ maps bijectively germs of monogenic sections at $x'$ to germs of monogenic sections at $x$. So it is enough to understand the germs of  monogenic sections at a fixed point  $x_0$ (which we will call the \textit{origin}). We will choose an open affine neighbourhood  $\aff=A(k,\R)\times M(n+1,k,\R)$ of $x_0$ where $A(k,\R)$ is the space of skew-symmetric matrices of size $k\times k$.  The set $\aff$ inherits in a non-canonical way the structure of a Lie group with a 2-graded nilpotent Lie algebra. The naturality of $D$ means that this operator is left invariant and so it commutes with right invariant vector fields which are the infinitesimal symmetries. Fixing trivializations of natural vector bundles over $\aff$, the $k$-Dirac operator can be viewed as a differential operator 
\begin{equation}
D:\cC^\infty(\aff,\Sp)\ra \cC^\infty(\aff,\R^k\otimes\Sp)
\end{equation}
with polynomial coefficients (see the formula (\ref{Dirac operator locally})).  From the formula (\ref{Dirac operator locally}) also immediately follows that there is a commutative diagram
\begin{equation}
\xymatrix{\cC^\infty(\aff,\Sp)\ar[r]^D& \cC^\infty(\aff,\R^k\otimes\Sp)\\
\cC^\infty(\aff,\Sp)\ar[u]^{\rho^\ast}\ar[r]^E&\cC^\infty(\aff,\R^k\otimes\Sp)\ar[u]^{\rho^\ast}\\}
\end{equation}
where we have for a moment denoted by $\rho$ the canonical projection $\aff\ra M(n+1,k,\R)$. So $D$ can be viewed as  a non-trivial "extension" of $E$ to the larger affine set $\aff$ and this also  justifies the name for $D$. As $\rho^\ast$ is injective, we know that  $\psi$ is a monogenic function iff  $D(\rho^\ast\psi)=0$.

As $D$ is a differential operator with polynomial coefficients, it is no longer true that homogeneous components of a monogenic spinor are again monogenic spinors.  However, this can be fixed by introducing (see Section \ref{section affine subset}) the weighted degree of polynomials. We will  write $\wdeg(f)=i$ if $f$ is a vector valued on $\aff$ such that each component of $f$ is a polynomial which is homogeneous of the weighted degree $i$ and put $\ms_i:=\{\Psi\in\cC^\infty(\aff,\Sp):D\Psi=0,\wdeg(\Psi)=i\}$. Then $D$ is a homogeneous operator of degree $-1$, i.e. $\wdeg(D(\Psi))=i-1$ whenever  $\wdeg(\Psi)=i$. It follows that $\Psi$ is a real analytic monogenic spinor iff $\Psi$ is a convering sum $\sum_{i\ge0}\Psi_i$ with $\Psi_i\in\ms_i$.

If we fix the origin and apply the classical Cartan-K\"ahler theorem to $D$, then we would actually view $D$ as the constant coefficient differential operator $E$ with its trivial extension to the larger affine set $\aff$ and this turns out not to be a good move (see \cite{S}). Rather than this, one could try to adapt the Cartan's lemma to the concept of homogeneous functions with respect to the weighted degree. In particular, we need to define a filtration on the spaces $\ms_i$ (more precisely on spaces which are isomorphic to $\ms_i$ which have more invariant meaning). Taking into account that the right invariant vector fields on $\aff$ are infinitesimal symmetries of $D$, it is clear that the filtrations should be defined with respect to these vector fields. However, there is a problem that the right invariant vector fields do not commute and so in general differentiation by a right invariant vector field does not give map $\ms_i\ra\ms_{i-1}$ that is compatible with the filtrations. Nevertheless, we will show that there is an ordered basis (see Lemma \ref{lemma preferred basis}) of the vector space of right invariant vector fields for which we get maps that are compatible with the filtrations. Here we will need some particular properties of the operator $D$  (this is used in the proof of Lemma \ref{lemma help}) and so at this point the machinery presented in this paper fails to work for general left invariant differential operators without further assumptions.
With this at hand  we can imitate the proof of Proposition 2.5 from \cite{Br} that a prolongation of an involutive system is again involutive and this is the hardest step in proving the main result of the paper.

\begin{thm}\label{main thm}
Let $i$ be a non-negative integer  and $\psi$ be a $\Sp$-valued function defined on the subset $M(n,k,\R):=\{x_{n+1,1}=\ldots=x_{n+1,k}=y_{12}\ldots=y_{k-1,k}=0\}$ of $\aff$ such that each component of $\psi$ is  a homogeneous polynomial of the degree $i$. Then there is a unique monogenic spinor $\Psi\in\ms_i$ such that $\Psi|_{M(n,k,\R)}=\psi$.
\end{thm}
In other words, we will consider here only formal solutions of the $k$-Dirac operator and we will not go into analysis of real analytic or smooth solutions of this operator.

Replacing the degree of polynomials by the weighted degree corresponds in more invariant language to replacing usual jets by weighted jets. Similarly, one needs to replace the  order of a differential operator by the weighted order. These concepts (see Section \ref{section weighted jets}) were developed in the 90's by Morimoto in \cite{MoI}, \cite{MoII} and \cite{MoIII}. See also \cite{K} and \cite{Ne} for an  application of  this  theory.  Many parts of the classical theory of exterior differential systems and the Cartan-K\"ahler theorem was generalized to the setting of weighted differential operators  in  \cite{MoIII}, \cite{MoIV} and \cite{MoV}. Namely, we will need the Spencer complex for weighted differential operators which has been introduced already in \cite{Ta}.

I would like to thank Katharina Neusser and the unknown referee for useful comments and to  Boris Krukligov and Peter Vassiliou for their kind support.

\subsection{Notation}

Let $p:V\ra M$ be a fibre  bundle. We denote the fibre  $p^{-1}(x)$ over $x\in M$ by $V_x$. If $U\subset M$ is an open subset, we denote  $p^{-1}(U)$ also by $V|_U$.
\bigskip

\textbf{Used symbols}.\\
$I(n,k)=\{(\alpha,i):\alpha=1,\dots,n,i=1,\dots,k\}$\\
$\hat I(k)=\{(r,s):r,s\in\Z;1\le r<s\le k\}$\\
$M(n,k,\R)$\ matrices of size $n\times k$ with real coefficients\\
$A(k,\R)$\ skew-symmetric matrices of size $k\times k$ with real coefficients\\
$1_k$\ is the indentity $k\times k$ matrix\\
$[v_1,\dots,v_\ell]$\ the linear span of vectors $v_1,\dots,v_\ell$\\
$\mV^\ast$, resp. $V^\ast$\ the dual of vector space $\mV$, resp. of vector bundle $V$\\

\section{The parabolic geometry}\label{section geometric structure}
In this section we will recall some basic knowledge from the theory of parabolic geometries and talk about the geometric structure for which the $k$-Dirac operator is natural. In the second part of the section we will go to an affine subset of the homogeneous model of the geometry, we will set some notation and definitions. Most importantly, we will introduce the weighted degree of polynomials on this affine set.

Let $\{e_1,\dots,e_k,\varepsilon_1,\dots,\varepsilon_{n+1},e^1,\dots,e^k\}$ be the standard basis of $\R^{2k+n+1}$. Then the bilinear form $h$ on $\R^{2k+n+1}$ which satisfies  $$h(e_i,e^j)=\delta^j_i,\ h(\varepsilon_\alpha,\ \varepsilon_\beta)=\delta_{\alpha\beta},\ h(e_i,e_j)=h(e^i,e^j)=0$$ where  $\delta$ is the Kronecker delta and $i,j=1,\dots,k;\alpha,\beta=1,\dots,n+1$, is non-degenerate and the associated quadratic form $H$ has signature $(k,k+n+1)$. We will sometimes write $\R^{k,n+k+1}$ instead of $\R^{2k+n+1}$ to indicate that the vector space comes with the bilinear form $h$. 

Let $\hat\G$ be  the   associated special orthogonal group $\SO(k,n+k+1)$.  The Lie algebra of $\hat\G$ is the simple  matrix algebra

\begin{equation}\label{orthogonal Lie alg}
\lag=\Bigg\{
\left(\begin{array}{ccc}
A&Z^T&W\\
X&B&-Z\\
Y&-X^T&-A^T
\end{array}\right)\Bigg|
\begin{matrix}
A\in M(k,\R),B\in A(n+1,\R),\\
X,Z\in M(n+1,k,\R),Y,W\in A(k,\R)
\end{matrix}
\Bigg\}.
\end{equation}
The block decomposition determines the direct sum decomposition $\lag_{-2}\oplus\lag_{-1}\oplus\lag_0\oplus\lag_1\oplus\lag_2$ of $\lag$ which is given by   

\begin{equation}\label{real lie algebra g}
\left(
\begin{array}{ccc}
\lag_0&\lag_1&\lag_2\\
\lag_{-1}&\lag_0&\lag_1\\
\lag_{-2}&\lag_{-1}&\lag_0
\end{array}
\right).
\footnote{Here we mean the following: $\lag_0$ is the subspace of block diagonal matrices, $\lag_1$ lives in the blocks $(1,2)$ and $(2,1)$   etc.}
\end{equation} 
If $X\in\lag_i,Y\in\lag_j$, then $[X,Y]\in\lag_{i+j}$ where we agree that $\lag_s=\{0\}$ if $|s|\ne0,1,2$. Moreover, $\lag_{-1}$ generates $\lag_-:=\lag_{-2}\oplus\lag_{-1}$ as Lie algebra. This means  that the decomposition is a $|2|$-gradation (see Definition 3.1.2 in \cite{CS}) on  $\lag$. The associated filtration of $\lag$ is $\{\lag^i:=\oplus_{j\ge i}\lag_j;i\in\Z\}$. Then  $\lap:=\lag^0$ is  a subalgebra and  each subspace $\lag^i$  is $\lap$-invariant.  Moreover $\lap=\lag_0\oplus\lap_+$ where  $\lap_+:=\lag^1$ is a nilradical  (i.e. a maximal nilpotent subalgebra) and  $\lag_0$ is a reductive Levi factor  of $\lap$ (see Section 2.1.8 in \cite{CS}). Notice that $\lag_0$ is isomorphic to  $\lagl(k,\R)\oplus\laso(n+1)$ and that the subspaces $\mE:=[e_1,\dots,e_k],\mF:=[\varepsilon_1,\dots,\varepsilon_{n+1}]$ are $\lag_0$-invariant. The space $\mF$ has $\lag_0$-invariant inner product  $h|_\mF$ which yields isomorphism $\mF\cong\mF^\ast$ as  $\lag_0$-module.  We  have the following list of isomorphisms of $\lag_0$-modules

\begin{equation}\label{g0 pieces in lap}
\lag_{-1}\cong\mE^\ast\otimes\mF,\ \lag_{-2}\cong\Lambda^2\mE^\ast\otimes\R,\ \lag_1\cong\mE\otimes\mF,\ \lag_2\cong\Lambda^2\mE\otimes\R
\end{equation}
where $\R$ is a trivial representation of $\laso(n)$. 

Let us now look for a homogeneous model of the geometry. The  group $\hat\G$ has a canonical transitive action on the Grassmannian variety $M$ of $k$-dimensional totally isotropic subspaces in $\R^{k,n+k+1}$. We  denote the action by dot, i.e. if $x\in M,g\in\hat\G$ then $g.x\in M$. In particular, the subspace $\mE$ is totally isotropic and we will call this point the \textit{origin of} $M$ and denote it by $x_0$. The stabilizer of $x_0$  is  a closed subgroup of $\hat\G$ with Lie algebra $\lap$. This group, we denote it by $\hat\gP$, is called a parabolic subgroup of $\G$ corresponding to the $|2|$-gradation.   The group  $\hat\gP$ is (see Proposition 3.1.3 in \cite{CS}) isomorphic to a semidirect product $\hat\G_0\ltimes\exp(\lap_+)$ with normal subgroup $\exp(\lap_+)$ where  $\hat\G_0$ is the subgroup of $\hat\G$ of the block diagonal matrices. Notice that $\hat\G_0$ is isomorphic to $\GL(k,\R)\times\SO(n+1)$ and that its Lie algebra is obviously  $\lag_0$.  So any $\hat\G_0$-module is naturally also a $\hat\gP$-module with  trivial action of $\exp(\lap_+)$ and it is well known (see Proposition 3.1.12 from \cite{CS}) that  any irreducible $\hat\gP$-module is of this form. This in particular applies to the nilpotent graded Lie algebra $\lag_-$ which becomes  a (non-irreducible) $\hat\gP$-module. Then the canonical map $\lag_-\ra gr(\lag/\lap)$ where $gr(\lag/\lap)$ is the  graded vector space associated to the  filtration $\lag^{-1}/\lap\subset\lag/\lap$ is an isomorphism of $\gP$-modules. 

The canonical projection  $\hat p:\hat\G\ra M,\hat p(g):=g.x_0$ is a principal $\hat\gP$-bundle over $M$. 
The quotient $\hat\cG_0$ of $\hat\G$ under the principal action of $\exp(\lap_+)$ is (see Section 3.1.5 in \cite{CS}) the total space of the principal bundle $\hat p_0:\hat \cG_0\ra M$ with  typical fibre $\hat\G_0\cong\hat\gP/\exp(\lap_+)$ where $\hat p_0(g.\exp(\lap_+)):=\hat p(g),g\in\G$. It is a well known fact that the tangent bundle $TM$ is isomorphic (see  Section 1.4.3  in \cite{CS}) to  the associated vector bundle  $\hat\G\times_{\hat\gP}(\lag/\lap)$.   This bundle contains the subbundle $\hat\G\times_{\hat\gP}(\lag^1/\lap)$ and so we see that $TM$ contains a canonical subbundle which we denote by $H$. Put $Q:=TM/H$ and let $gr(TM)$ be the associated graded vector bundle.  The bundle $gr(TM)$ is  isomorphic to $\hat\G\times_{\hat\gP} gr(\lag/\lap)$. We will view $gr(\lag/\lap)$ rather as $\lag_-$. As $\exp(\lap_+)$ acts trivially on $\lag_-$, it follows that $gr(TM)$ is isomorphic to $\hat\cG_0\times_{\hat\G_0}\lag_-$. We see that for each $x\in M$ the graded vector space $gr(T_xM)$ has the structure of a nilpotent graded 
Lie algebra. This structure  has a geometric origin, i.e. it is  well known that the structure   coincides with the Levi bracket  $\cL_x:\Lambda^2 gr(T_xM)\ra gr(T_xM)$ which is naturally induced by the Lie bracket of vector fields (see Section 3.1.7 in \cite{CS}). In particular, we see that $H$ is a bracket generating distribution.

The bundle $\hat\cG_0$ can be viewed also as the frame bundle of  a tautological vector bundle over $M$.
Let $E$, resp. $E^\bot$ be the total space of the  vector bundle over $M$ such that $E_x=x$, resp. $E^\bot_x=x^\bot$. Then $E\subset E^\bot\subset M\times\R^{2k+n+1}$ and the vector bundle $F:=E^\bot/E$ has a canonical inner product which comes from restricting $h$ to its fibres.  The quotient bundle $(M\times\R^{2k+n+1})/E^\bot$ is canonically isomorphic to $E^\ast$ where the canonical pairing with $E$  is again naturally induced by $h$. As the canonical volume form on $\R^{2k+n+1}$ is invariant under $\hat\G$ and the bundle $E\oplus E^\ast$ has a canonical orientation, we can fix the orientation on the bundle  $F$  so that for each $x\in M$ the induced orientation on $E_x\oplus F_x\oplus E^\ast_x$ coincides with the canonical one on $\R^{2k+n+1}$. Obviously $E$, resp. $F$ is a $\hat\G$-homogeneous vector bundle whose fibre  over the origin $x_0$ is isomorphic to the (irreducible) $\gP$-module  $\mE$, resp. $\mF$. It follows (see Section 1.4.3 in \cite{CS}) that $E$, resp. $F$ is  isomorphic to the associated vector bundle $\hat\G\times_{\hat\gP}\mE$, resp. $\hat\G\times_{\hat\gP}\mF$ and since $\exp(\lap_+)$ acts trivially on  $\mE$, resp. $\mF$,  it follows that it is  isomorphic to $\hat\cG_0\times_{\hat\G_0}\mE$, resp. $ \hat\cG_0\times_{\hat\G_0}\mF$. Hence, we may view $(\hat\cG_0)_x,x\in M$ as the set of pairs $(p,q)$ where $p$ is a frame of $E_x$ and $q$ is an orthonormal frame of $F_x$ which is compatible with the orientation on $F_x$. Repeating this argument also to the bundles $H,Q,H^\ast,Q^\ast$ respectively and using (\ref{g0 pieces in lap}), we obtain isomorphisms
\begin{equation}\label{associated vector bundles}
H\cong E^\ast\otimes F,\ Q\cong\Lambda^2 E^\ast,\ H^\ast\cong E\otimes F,\ Q^\ast\cong\Lambda^2E.
\end{equation}

 The group $\hat\G$ is not simply connected. In order to invariantly define the $k$-Dirac operator we will choose a 4:1 covering  $\rho:\G\ra\hat\G$. It is easy to check that the inclusion $\hat\G_0\hookrightarrow\G$ induces isomorphisms on $\pi_0$  and  $\pi_1$.
It follows that 
\begin{eqnarray}
\pi_1(\hat\G)=
\bigg\{
\begin{matrix}
\ \Z\times\Z_2,\ k=2,\\
\ \Z_2\times\Z_2,\ k>2.
\end{matrix}
\end{eqnarray}
If $k=2$, then $\G$ is determined by the subgroup  $2\Z\times\{1\}$. If $k>2$, then $\G$ is a universal covering of $\hat\G$. We put $\gP:=\rho^{-1}(\hat\gP),\G_0:=\rho^{-1}(\hat\G_0)$. It follows that  $\G_0$ is isomorphic to $\widetilde\GL(k,\R)\times\Spin(n+1)$ where $\widetilde\GL(k,\R)$ is isomorphic to a connected  2:1 covering of $\GL(k,\R)$ (which is unique as the fundamental group of $\GL(k,\R)$ contains a unique subgroup of index $2$). We obtain a principal fibre bundle $p:\G\ra M$ with typical fibre $\gP$ where we put $p:=\hat p\circ\rho$. We define $\cG_0:=\G/\exp(\lap_+)$ so that the map $\rho$ factorizes to a 4:1 covering  $\rho_0:\cG_0\ra\hat\cG_0$.  We get another  principal  bundle $p_0:\cG_0\ra M$ where $p_0:=\hat p_0\circ\rho_0$ with typical fibre $\G_0$.

 Let $\Sp$ be the  space of spinors for the complexified Clifford algebra $(\mF,h|_{\mF})$. We may view the pair $(\mF,h|_\mF)$ as well  as the pair $(\R^{n+1},g)$ from the introduction and so $\Sp$ is the same of spinors. Recall that if $n+1$ is odd, then $\Sp$ is an irreducible representation of $\Spin(n+1)$ while if $n+1$ is even, then $\Sp$ is  the direct sum of two irreducible subspaces. Let $\C_\lambda$ be an irreducible complex $\GL(k,\C)$-module with highest weight $\lambda:=(\frac{n-1}{2},\frac{n-1}{2},\ldots,\frac{n-1}{2})$ as in \cite{Sa}. Then $\C_\lambda$ is also a 1-dimensional $\GL(k,\R)$-module by restriction on which only the center (the subspace of multiples of the identity matrix) acts non-trivially.  We get two irreducible $\G_0$-modules  $\Sp_\lambda:=\C_\lambda\otimes_\C\Sp$ and $\Sp^k_\lambda:=(\mE\otimes_\R\C_\lambda)\otimes_\C\Sp$ where we indicate by subscript over which field we take the tensor product. These two space are also irreducible $\gP$-modules as we explained above. If we forget the structure of $\G_0$-modules 
on both spaces, then there are canonical isomorphisms 
\begin{equation}\label{isomorphism of spinor modules}
\Sp_\lambda\ra\Sp,\ \mathrm{resp.} \ \ \Sp^k_\lambda\ra\mE\otimes_\R\Sp
\end{equation} 
of  vector spaces that are given by $z\otimes\psi\mapsto z.\psi$, resp. $(v\otimes z)\otimes\psi\mapsto v\otimes (z.\psi)$ where $z\in\C_\lambda,v\in\mE,\psi\in\Sp$.  Now we can form two associated vector bundles:
 \begin{equation}\label{spinor bundles}
 S_\lambda:=\G\times_\gP\Sp_\lambda,\ S^k_\lambda:=\G\times_\gP\Sp^k_\lambda.
 \end{equation}
The $k$-Dirac operator then maps section of $S_\lambda$ to sections of $S^k_\lambda$. We will give an invariant definition of $D$ in 
 Section \ref{section k Dirac operator}. Later on, we will trivialize the bundles  over the affine subset $\aff$ and we will give a formula $D$ with respect to these trivializations.

\subsection{The affine set $\aff$}\label{section affine subset}
Recall that 
\begin{equation}\label{lie algebra g_-}
\lag_{-}=\Bigg\{
\left(
\begin{matrix}
0&0&0\\
X&0&0\\
Y&-X^T&0\\
\end{matrix}
\right)\Bigg|
X=\in M(n+1,k,\R),Y=-Y^T\in A(k,\R)
\Bigg\}
\end{equation}
and that $\lag_-=\lag_{-2}\oplus\lag_{-1}$ where  $\lag_{-2}$, respectively $\lag_{-1}$ is the subspace of those matrices where  $X=0$, resp. $Y=0$.   The map  $p\circ\exp:\lag_-\ra M$ is
\begin{eqnarray}\label{coordinates on osu}
&p\circ\exp\left(
\begin{array}{ccc}
0&0&0\\
X&0&0\\
Y&-X^T&0\\
\end{array}
\right) 
=
\left[
\begin{matrix}
1_k\\
X\\
Y-\frac{1}{2}X^TX\\
\end{matrix}
\right]&
\end{eqnarray}
where $[a_{ij}]$ is the $k$-plane in $\R^{2k+n+1}$ spanned by the columns of a matrix $(a_{ij})\in M(2k+n+1,k,\R)$. We see that the map  (\ref{coordinates on osu}) is  injective  and that its image $\aff$ is an affine subset of $M$. It is also easy to see that  $\G_-:=\exp(\lag_-)$ is a closed connected analytic subgroup of $\G$ with Lie algebra $\lag_-$ such that the exponential map $\lag_-\ra\G_-$ is bijective.   Let us write the matrices above as $X=(x_{\alpha i})_{\alpha=1,\dots,n+1}^{i=1,\dots,k}$, resp. $Y=(y_{rs})_{r=1,\dots,k}^{s=1,\dots,k}$. We will view  $\lag_-$ also as an affine space with coordinates
\begin{equation}\label{coordinates}
x_{\alpha i},y_{rs}\ \mathrm{where}\ (\alpha, i)\in I(n+1,k),(r,s)\in \hat I(k).
\end{equation}
 
\begin{remark}
 We see that both maps in the composition 
 \begin{equation}\label{isomorphisms}
\lag_-\xrightarrow{\exp}\G_-\xrightarrow{p|_{\G_-}}\aff 
 \end{equation}
are diffeomorphisms which allows us to push-forward and pullback geometric objects from one set to another set. We will do that without any further comment. In particular, we will use the matrix coefficients on $\lag_-$  also as coordinates on  $\aff$ and we will view left and right invariant vector fields on $\G_-$ as vector fields on $\aff$. 
\end{remark}

Let $X\in\lag_-$. We denote by $L_X$, resp. $R_X$ the corresponding left, resp. right invariant vector field on $\G_-$. As the distribution $H$ is invariant with respect to the canonical left action of $\G$ on $M$, it follows that $H$ is spanned over $\aff$ by the vector fields $L_X$ as $X$ ranges over the set $\lag_{-1}$. It is a straightforward calculation  to show that
\begin{eqnarray}\label{left invariant vector fields}
&&L_{\alpha i}:=L_{e_i\otimes\varepsilon_\alpha}=\partial_{x_{\alpha i}}-\frac{1}{2}\sum_{j=1}^kx_{\alpha j}\partial_{y_{ij}}, \\ 
&&R_{\alpha i}:=R_{e_i\otimes\varepsilon_\alpha}=\partial_{x_{\alpha i}}+\frac{1}{2}\sum_{j=1}^kx_{\alpha j}\partial_{y_{ij}}\label{right invariant vector fields},\\
&&L_{e_i\wedge e_j}=R_{e_i\wedge e_j}=\partial_{y_{ij}}
\end{eqnarray} 
where we  use the isomorphisms from (\ref{g0 pieces in lap}) and the bases of $\mE,\mF$ that are given above that formula. Moreover, we will use the convention $\partial_{y_{rs}}=-\partial_{y_{sr}}$. We have that
\begin{equation}\label{Lie bracket}
[L_{\alpha i},L_{\beta j}]=-[R_{\alpha i},R_{\beta j}]=\delta_{\alpha\beta}\partial_{y_{ij}}
\end{equation}
while all  other Lie brackets of the vector fields given above are zero. 

The left invariant vector fields form a natural framing over $\aff$ which is adapted to the filtration $H\subset TM$.

 \begin{ex}$($\textsl{Left invariant framing and coframing over $\aff$}$)$ \label{example coframing}
 
We see that  $\{L_{\alpha i},\partial_{y_{rs}}:(\alpha,i)\in I(n+1,k),(r,s)\in \hat I(k)\}$ is a framing over $\aff$ which is adapted to the filtration of $T\aff$ as for each $(\alpha,i)\in I(n+1,k)$ the vector field $L_{\alpha i}$ is a section of $H|_\aff$.  Let $\{\omega_{\alpha i},\theta_{rs}\}$ be the dual coframing, i.e.  $\omega_{\alpha i},\theta_{rs}$ are differential 1-forms on $\aff$  such  that 
\begin{eqnarray*}
&\omega_{\beta j}(L_{\alpha i})=\delta_{\alpha\beta}\delta_{ij},\ \theta_{rs}(\partial_{y_{ij}})=\delta_{ri}\delta_{sj}-\delta_{rj}\delta_{si},\ \omega_{\alpha i}(\partial_{y_{rs}})=\theta_{rs}(L_{\alpha i})=0\nonumber&
\end{eqnarray*}
 Then we have 
\begin{equation}
\omega_{\alpha i}=dx_{\alpha i},\ \theta_{rs}=dy_{rs}-\frac{1}{2}\sum_{\alpha=1}^{n+1}(x_{\alpha r}dx_{\alpha s}-x_{\alpha s}dx_{\alpha r}).
\end{equation} 
In particular, the exterior derivative of $\theta_{rs}$ is
\begin{equation}\label{ex of left inv contact form}
d\theta_{rs}=-\sum_{\alpha=1}^{n+1}\omega_{\alpha r}\wedge\omega_{\alpha s}.
\end{equation}
 \end{ex}

The framing over $\aff$ determined by the right invariant vector fields is not adapted to the filtration of the tangent bundle of $M$. On the other hand, the right invariant vector fields are infinitesimal symmetries of the parabolic structure and this is  why we will use these vector fields later on.

\begin{ex}$($\textsl{Right invariant framing and coframing over $\aff$}$)$

Let $\{\varpi_{\alpha i},\vartheta_{rs}\}$ be the dual coframing to the framing $\{R_{\alpha i},\partial_{y_{rs}}\}$ over $\aff$  where $(\alpha, i),(r,s)$ ranges over the sets $ I(n+1,k),\hat I(k)$ respectively. 
We find that 
\begin{eqnarray}
\varpi_{\alpha i}=dx_{\alpha i},\ \vartheta_{rs}=dy_{rs}+\frac{1}{2}\sum_{\alpha=1}^{n+1}(x_{\alpha r}dx_{\alpha s}-x_{\alpha s}dx_{\alpha r})
\end{eqnarray}
and in particular
\begin{equation}\label{differential of contact forms}
d\vartheta_{rs}=\sum_{\alpha=1}^{n+1}\omega_{\alpha r}\wedge\omega_{\alpha s}
\end{equation}



\end{ex}

The $k$-Dirac operator $D$ is not a constant coefficient  operator but it is  a differential operator with   polynomial coefficients (see the formula (\ref{Dirac operator locally}) but have in mind that this formula depends on  choices). Hence,  it is not true that homogeneous components of a real analytic monogenic spinor  are again monogenic spinors. However, we may fix this by introducing a weighted degree. 

Let $\R[x,y]$, resp. $\R[y]$, resp. $\R[x]$ be the ring of polynomials with real coefficients on the affine space $\lag_-$, resp. $\lag_{-2}$, resp. $\lag_{-1}$. Then since $\lag_{-2}$, resp. $\lag_{-1}$ is an  affine subspace of $\lag_-$ we may view $\R[y]$, resp. $\R[x]$ as a subring of $\R[x,y]$. If $s\in\R[x,y]$ is a monomial, then we can certainly find monomials $p'\in\R[y],p''\in\R[x]$ such that $p'.p''=s$. Then the number $\wdeg(s):=2\dg(p')+dg(p'')$  where $\dg$ is the usual degree of polynomials is independent of the choice of $p',p''$ and we call it the  \textit{weighted degree of} s. More generally, we call a polynomial $p\in\R[x,y]$ a \textit{homogeneous polynomials of the weighted degree} $i$ (and  write $\wdeg(p)=i$) if $p$ is a sum of monomials of the weighted degree $i$. We denote the vector space of all homogeneous polynomials of the weighted degree $i$ by $\R[x,y]_i$ and put $\R[x,y]_i=\{0\}$ if $i<0$. Then we clearly have $\R[x,y]=\oplus_{i\ge0}\R[x,y]_i$. Also notice that the multiplication is compatible with 
the grading. 

\begin{ex}\label{example homogeneous polynomials}
 $y_{12}+x_{11}$ is not a homogeneous polynomial with respect to the weighted degree as $\wdeg(y_{12})=2$ and $\wdeg(x_{11})=1$. Also notice that if $p$ is a homogeneous polynomial of the weighted degree $j$,  then $L_{\alpha i}p,R_{\alpha i}p\in\R[x,y]_{j-1}$ and $\partial_{y_{rs}}p\in\R[x,y]_{j-2}$.
\end{ex}

Let $\mU$ be a vector space. We say that a $\mU$-valued function $f$ on $\cA$ is a \textit{homogeneous function of a weighted degree} $i$ (and write $\wdeg(f)=i$) if each component of $f$ (with respect to a fixed basis of $\mU$) belongs to $\R[x,y]_i$. Notice that the definition does not depend on the choice of the basis of $\mU$. 

\section{Weighted jets and weighted differential operators on filtered manifolds}\label{section weighted jets}

Now we will recall the notion of weighted differential operators and weighted jets (see also \cite{K}, \cite{MoI}, \cite{MoII}, \cite{MoIII} and \cite{Ne} for a broader introduction to the concept). We will see that this a natural concept for filtered  manifolds which generalizes the notion of the usual degree of differential operators and the jets of the germs of sections of vector bundles. We will also see that the weighted jet of a function at a point is an invariant analogue of the weighted degree of functions that we defined in the previous section (see Example \ref{ex graded jets} below that the map which assigns to a function its weighted jet at the origin of $M$ gives an isomorphism between the space of polynomials graded with respect to the weighted degree and the graded space of weighted jets of functions.) In the second part of the section we will give an invariant definition of the $k$-Dirac operator and we will recall the definition of the tableau and its prolongations (see also material \cite{MoIII} and \cite{Ne} where this can be found).

Recall that there is a natural bracket generating distribution $H$ on $M$ such that  $(gr(TM),\cL)$ is a locally trivial bundle of graded nilpotent Lie algebras over $M$ with  a typical fibre $\lag_-$. Put $F_{-1}:=H,F_{-2}:=TM$. Let $X$ be a vector field which is defined over an open subset $U$ of $M$, then the \textit{weighted order of} $X$  is defined as the smallest integer $i$ such that $X\in\Gamma(F_{-i}|_U)$. We write $ord(X)=i$.

\begin{ex}\label{example graded jets of functions at origin}
We have that $ord(L_{\alpha i})=1,ord (R_{\alpha i})=ord(\partial_{y_{rs}})=2$ where the vector fields were defined in (\ref{left invariant vector fields}) and $(\alpha,i)\in I(n+1,k),(r,s)\in \hat I(k)$.   
\end{ex}

 A differential operator $D$ acting on the space of smooth functions on $M$ is called a \textit{differential operator of the weighted order at most} $r$  if for each point $x\in M$ there is an open neighbourhood $U$ of $x$ with  local framing $\{X_1,\dots,X_p\}$ such that 
\begin{equation}\label{weighted differential operator}
D|_U=\sum_{a\in\mathbb N_{0}^l} f_a X_1^{a_1}\dots X_p^{a_p}
\end{equation}
where $\mathbb N_0^l:=\{a=[a_1,\dots,a_p]:a_i\in\Z,a_i\ge0,i=1,\dots,p\},f_a\in\cC^\infty(U)$ and for all $a$ in the sum with $f_a$ non-zero: $\sum_{i=1}^pa_i.ord(X_i)\le r$. If $D$ is a differential operator of the order at most $r$ but not of the order at most $r-1$ then we say that $D$ is a \textit{differential operator of the weighted order} $r$ (and we will write $ord(D)=r$). 
 
Let $f,f'$ be two germs of  smooth functions at $x\in M$. Then we say that $f,f'$ are $r$-equivalent $($we write $f\sim_r f')$ if $Df(x)=Df'(x)$ for all differential operators of the weighted order at most $r$. Clearly,  $\sim_r$ is an equivalence relation and we denote by $\xj^r_xf$ the equivalence class of the germ $f$ and by $\xJ^r_x$ the space of all such equivalence classes. The disjoint union $\xJ^r:=\cup_{x\in M}\xJ_x^r$ is naturally a vector bundle over $M$. There is a canonical projection $\pi_r:\xJ^r\ra\xJ^{r-1}$  whose kernel $\mathfrak{gr}^r$ is again a vector bundle.

Let $V$ be a vector bundle over $M$ and  $s,s'$ be two germs of smooth sections of $V$ at  $x\in M$. We say that $s,s'$ are $r$-equivalent ($s\sim_r s$) if
\begin{equation}
D\langle\lambda,s-s'\rangle(x)=0
\end{equation}
for all sections $\lambda$ of the dual bundle $V^\ast$  and all differential operators $D$ of the weighted order at most $r$. Here $\langle-,-\rangle$ denotes the canonical pairing. We again denote by $\xj^r_xs$  the equivalence class of the germ $s$ and  the space of all such equivalence classes by $\xJ^r_xV$. The disjoint union $\xJ^rV:=\bigcup_{x\in M}\xJ^r_xM$ is  a vector bundle over $M$ (see for example Theorem 2.10 in \cite{Ne}) such that the canonical map $\pi_r:\xJ^rV\ra\xJ^{r-1}V$ is clearly a surjective vector bundle map whose kernel $\gr^rV$ is again a vector bundle over $M$. We denote its fibre over $x\in M$ by $\gr^r_x V$.

Let us now assume that $V$ is a $\G$-homogeneous vector bundle. Then $V$ is naturally isomorphic to $\G\times_\gP\mV$ where  $\mV$ is the $\gP$-module $V_{x_0}$ (see Section 1.4.3 in \cite{CS}). The bundles   $\xJ^rV,\gr^rV$ have also a canonical $\G$-action and so  it suffices to understand their fibres (see the following example) over the origin $x_0$.   We will use the following notation:
\begin{equation}\label{notation for fibres of weighted jets over the origin for homogeneous bundles}
\xJ^r\mV:=\xJ^r_{x_0}\mV,\gr^r\mV:=\gr^r_{x_0}\mV.
\end{equation}

\begin{ex}\label{ex graded jets}
Recall that $\lag_-=\lag_{-2}\oplus\lag_{-1}$ is a nilpotent graded Lie algebra, i.e. $[\lag_{i},\lag_{j}]\subset\lag_{i+j}$ where $i,j=-1,-2$ and we put $\lag_{\ell}=\{0\}$ whenever $\ell\ne-1,-2$. The universal enveloping algebra $\cU(\lag_-)$ of $\lag_-$ is defined as $T(\lag_-)/I$ where $T(\lag_-)$ is the tensor algebra of $\lag_-$ and $I$ is the both sided ideal  which  is generated by the elements of the form $X\otimes Y-Y\otimes X-[X,Y];X,Y\in\lag_-$. As the Lie bracket on $\lag_-$ is compatible with the grading it follows that $I$ is spanned as a vector space by homogeneous elements. It follows that the enveloping algebra $\cU(\lag_-)$ is naturally graded  $\oplus_{r\ge0}\ \cU_{-r}(\lag_-)$ where $\cU_{-r}(\lag_-):=T_{-r}(\lag_-)/(T_{-r}(\lag_-)\cap I)$ and $T_{-r}(\lag_-)$ is the vector subspace of  $T(\lag_-)$ which is spanned by the elements of the form $X_{i_1}\otimes\dots \otimes X_{i_u}$ where  $X_{i_j}\in\lag_{s_j};s_j=-1,-2;j=1,\dots,u$ and  $\sum_{j=1}^us_j=-r$. 

Recall that we view $\lag_-$ as the graded space $gr(\lag/\lap)$  which is  canonically isomorphic to  $gr(T_{x_0}M)$.  We will denote the $j$-th graded component of $gr(T_{x_0}M)$ by $gr_j(T_{x_0}M)$ and so $gr_j(T_{x_0}M)$ is isomorphic to $\lag_j;j=-2,-1$. Given vectors $X_{i_1},\dots,X_{i_u}$  as above,  there are vector fields  $\hat X_{i_1},\dots,\hat X_{i_u}$  on $M$ such that  $ord(X_{i_j})=-s_j$ and 
\begin{equation}
X_{i_j}=
\bigg\{
\begin{matrix}
\hat X_{i_j}(x_0),s_j=-1,\\
q(\hat X_{i_j})(x_0),s_j=-2
\end{matrix}
\end{equation}
where $q:TM\ra Q$ is the canonical projection.  

Let $V$ be the $\G$-homogeneous vector bundle $\G\times_\gP\mV$. Then $V$ can be trivialized over the affine set  $\aff$  and so we can view the germ of a section $s$ of $V$ at $x_0$ as a vector valued function which we may differentiate by the vector fields $\hat X_{i_1},\dots,\hat X_{i_u}$. If $\xj^{r-1}_{x_0}s=0$, then the value
\begin{eqnarray}
&(\hat X_{i_1}\dots\hat X_{i_u} s)(x)&
\end{eqnarray}
depends only $\xj^r_{x}s$ and on $X_{i_1},\dots,X_{i_u}$ (in particular it does not depend on the  way we have extended the vector $X_{i_j}$ to the vector field $\hat X_{i_j},j=1,\dots,u$ and the choice of trivialization of $V$). Hence,  we have obtained a well defined map 

\begin{equation}\label{isomorphism of graded jets over origin}
\gr^r\mV\mapsto\cU_{-r}^\ast(\lag_-)\otimes\mV.
\end{equation}
The  construction above works for a general point $x\in M$ and we would get a map $\gr^r_xV\ra\cU_{-r}^\ast(gr(T_xM))\otimes V_x$ which is bijective by Proposition 2.2 from \cite{Ne}. 
Using the duality $\lag_{-i}^\ast\cong\lag_i$ then the isomorphism (\ref{isomorphism of graded jets over origin}) becomes for small $r$: 
$$\gr^1\mV\cong\lag_1\otimes\mV,\ \gr^2\mV\cong S^2\lag_1\otimes\mV\oplus\lag_2\otimes\mV,\dots$$
and in general
$$\gr^r\mV\cong\bigoplus_{i=0}^{\lfloor\frac{r}{2}\rfloor}S^{r-2i}\lag_1\otimes S^i\lag_2\otimes\mV$$
where $\lfloor a\rfloor $ is the integer part of $a\in\R$.

Finally, let $p\in\R[x,y]_r$.  Then (recall Example \ref{example homogeneous polynomials} and \ref{example graded jets of functions at origin}) we have that $\xj^{r-1}_{x_0}p=0$ and it is easy to verify that  the map $p\mapsto\xj^r_{x_0}p$ induces isomorphism $\R[x,y]_r\ra\gr^r_{x_0}$.  From (\ref{isomorphism of graded jets over origin}) we know that $\gr^r_{x_0}$ is isomorphic to $\cU^\ast_{-r}(\lag_-)$. Altogether, we  obtain isomorphism
\begin{equation}\label{weighted jets of spinors at the origin}
\R[x,y]_r\otimes\mV\ra\gr^r\mV.
\end{equation}
\end{ex}






Now we can give an invariant definition of the $k$-Dirac operator using the language of weighted jets.

\subsection{The $k$-Dirac operator and the tableau}\label{section k Dirac operator}

 There is (up to a constant) a  unique non-zero $\gP$-equivariant map
\begin{equation}	\label{P-eq map for k-Dirac operator}
\phi:\xJ^1\Sp_\lambda\ra\Sp^k_\lambda
\end{equation}
where $\Sp^k_\lambda,\Sp_\lambda$ where defined at the end of Section \ref{section geometric structure}. The map $\phi$ can be extended in a unique way to a $\G$-equivariant vector bundle map
\begin{equation}
\Phi:\xJ^1S_\lambda\ra S^k_\lambda.
\end{equation}
 For a non-negative integer $i$ we define the  $i$-th prolongation of $\Phi$ as the vector bundle map
\begin{eqnarray}
&\xp^i\Phi:\xJ^{1+i}S_\lambda\ra\xJ^iS^k_\lambda&\\
&\xp^i\Phi(\xj^{1+i}_xs):=\xj^i_x(\Phi(\xj^1s))&\nonumber
\end{eqnarray}
where $s$ is a germ of a section of $S_\lambda$ at $x$. As $\Phi$ is linear in fibres, the right hand side depends only the weighted  $(i+1)$-jet of $s$ and so $\xp^i\Phi$ is well defined (see also Section 1.3.1 in \cite{Ne}). We denote $\xp^i\phi$ the restriction of $\xp^i\Phi$ to the fibres over the origin and put 
\begin{equation}
\fs_i:=\bigg\{
\begin{matrix}
 \Sp_\lambda,\ \ i=0\\
 Ker(\xp^{i-1}\phi),\ i\ge 1.\\
\end{matrix} 
\end{equation} 
Then for each $i\ge 0$ there is a commutative diagram of vector spaces with exact rows and columns
\begin{equation}\label{big diagram}
\xymatrix{&0\ar[d]&0\ar[d]&0\ar[d]\\
0\ar[r]& \tab^{(i)}\ar[d]\ar[r]&\gr^{i+1}\Sp_\lambda\ar[r]\ar[d]&\gr^i\Sp^k_\lambda\ar[d]\\
0\ar[r]&\fs^{i+1}\ar[r]\ar[d]&\xJ^{i+1}\Sp_\lambda\ar[d]\ar[r]^{\xp^i\phi}&\xJ^i\Sp^k_\lambda\ar[d]\\
0\ar[r]&\fs^{i}\ar[r]&\xJ^{i}\Sp_\lambda\ar[r]^{\xp^{i-1}\phi}&\xJ^{i-1}\Sp^k_\lambda.\\}
\end{equation}
So  $\tab^{(i)}$ is at the same time the kernel of the canonical projection $\fs^{i+1}\ra\fs^{i}$ and  the kernel of the restriction of $\xp^i\phi$ to $\gr^{i+1}\Sp_\lambda$.  We will call $\tab:=\tab^{(0)}$ the \textit{tableau} (determined by $\phi$) and $\tab^{(i)}$ the $i$-\textit{th prolongation of the tableau} (see also \cite{MoIII}). It is shown in \cite{S} that 
\begin{eqnarray}
&&\tab=\mE_\lambda\otimes(\mF\boxtimes\Sp),\\
&&\tab^{(1)}=(S^2\mE\otimes\C_\lambda)\otimes (S^2_0\mF\boxtimes\Sp)\oplus(\Lambda^2\mE\otimes\C_\lambda)\otimes(\Lambda^2\mF\boxtimes\Sp\oplus\Sp)
\end{eqnarray}
where $\boxtimes$ is the highest weight component\footnote{In these special cases the highest weight components has a more direct description. It coincides with the kernel of $S^2_0\mF\otimes\Sp\ra\Sp$, resp. of $\Lambda^2\mF\otimes\Sp\ra\Sp$.} in the tensor product and the $S^2_0$ denotes the trace-free part of the second symmetric power. In particular, we easily compute that
\begin{equation}\label{dimension of tableaux}
\dim\tab=nk\dim\Sp_\lambda,\ \dim\tab^{(1)}={nk+1\choose 2}\dim\Sp_\lambda.
\end{equation}

 The $k$-Dirac operator can be then invariantly defined as the composition
\begin{equation}
\Gamma(S_\lambda)\ra\Gamma(\xJ^1S_\lambda)\xrightarrow{\Phi}\Gamma(S^k_\lambda)
\end{equation} 
where the first map is the canonical inclusion which assigns to a section of $S_\lambda$ its weighted $1$-jet at each point. As $\Phi$ is a $\G$-equivariant vector bundle map, it follows that $D$ is a $\G$-invariant linear differential operator. In particular, the operator $D$ commutes with any infinitesimal symmetry. Recall that infinitesimal symmetries are the right invariant vector fields. Also notice that   $Ds(x)=Ds'(x)$ whenever $\xj^1_xs=\xj^1_xs',x\in M$ and so we say that \textit{the weighted order of} $D$ is at most $1$. 

\begin{remark}\label{remark differentiating monogenic spinors}

Let $\hat X$ be a vector field of the weighted order $r$ which is defined on an open subset $U$ of $M$. Then the formula (\ref{isomorphism of graded jets over origin}) and the remark below the formula imply that differentiation by $\hat X$ induces a well defined map $\gr^i_xV\ra\gr^{i-r}_xV$ for each $x\in U,i\in\N$. In particular, for $X\in\lag_{-1}$ we know that $ord(L_X)=1$ and so  it follows that the differentiation by $L_X,X\in\lag_{-1}$ induces a linear map $\gr^i_xV\ra\gr^{i-1}_xV$ for each $x\in \aff$.
However, it is not true that if we choose $x=x_0$, that the map restricts to a map $\tab^{(i)}\ra\tab^{(i-1)}$. This follows from the fact that $L_X$ is not an infinitesimal symmetry of the structure and so the operator $D$ does not commute in general with $L_X$ (compare with the local formula (\ref{Dirac operator locally}) of the $k$-Dirac operator).

 On the other hand, by  Example \ref{example homogeneous polynomials} and (\ref{weighted jets of spinors at the origin}) it follows that  differentiation by the right invariant vector field $R_X$ induces a well defined map
\begin{equation}\label{differentiating with right invariant vector fields}
R_X:\gr^{i+1}\Sp_\lambda\ra\gr^i\Sp_\lambda.
\end{equation} 
  As the right invariant vector fields are infinitesimal symmetries of the parabolic structure, the map $R_X$ commutes with the $k$-Dirac operator $D$ (and hence, with all prolongations $\xp^i\phi$).    It follows that the map (\ref{differentiating with right invariant vector fields}) restricts to a well defined map
\begin{equation}
R_{X}:\tab^{(i)}\ra\tab^{(i-1)}.
\end{equation}
This observation will be crucial in the sequel. As the vector field  $\partial_{y_{rs}}$  is both  left and right invariant, it follows that differentiation by this vector fields induces a well defined linear map
\begin{equation}
\partial_{y_{rs}}:\tab^{(i)}\ra\tab^{(i-2)}
\end{equation}
Finally, notice that $s\in\gr^{i+1}\Sp_\lambda$ is contained in the subspace $\tab^{(i)},i\ge1$ iff $R_Xs\in\tab^{(i-1)}$ for each $X\in\lag_{-1}$.
\end{remark}

\section{The set of initial condition for the $k$-Dirac operator} \label{section k-Dirac operator}

Consider the map $\rho:\aff\xrightarrow{\mu}\G_-\hookrightarrow\G$ where $\mu$ is the inverse  of the map on the right hand side in (\ref{isomorphisms}) and the second map is the canonical inclusion. Then clearly $\rho$ is a smooth section of the canonical  $\gP$-bundle over $\aff$.
Given a section $\psi$ of the vector bundle $S_\lambda$  over $\aff$, there is  (see Proposition 1.2.7 in \cite{CS}) a unique $\gP$-equivariant function $f\in\cC^\infty(p^{-1}(\aff),\Sp_\lambda)^{\gP}$.  The map $\psi\mapsto \underline\psi:=f\circ\rho$ is a bijection $\Gamma(S_\lambda|_\aff)\ra\cC^\infty(\aff,\Sp_\lambda)$. As the structure of $\G_0$-modules is no longer visible here, we may view, using the isomorphism of the vector spaces from (\ref{isomorphism of spinor modules}), the latter space as $\cC^\infty(\aff,\Sp)$.   We similarly identify $\Gamma(\Sp^k_\lambda|_\aff)$ with $\cC^\infty(\aff,\R^k\otimes\Sp)$. Then there is a unique linear differential operator (which we also for simplicity also denote by $D$)
\begin{equation}\label{Dirac operator on functions}
D:\cC^\infty(\aff,\Sp)\ra\cC^\infty(\aff,\R^k\otimes\Sp)
\end{equation}
such that $D\underline\psi=\underline{D\psi}$ for each $\psi\in\Gamma(S_\lambda|_U)$. This is the $k$-Dirac operator in the trivialization induced by $\rho$. We will call $\Psi\in\cC^\infty(\aff,S)$ a \textit{spinor valued function} or simply a \textit{spinor}. Then (see \cite{Sa})
\begin{equation}\label{Dirac operator locally}
D\Psi=(D_1\Psi,\dots,D_k \Psi)\ \mathrm{where} \ D_i\Psi=\sum_{\alpha=1}^{n+1}\varepsilon_\alpha. L_{\alpha i}\Psi 
\end{equation}
and the dot denotes the Clifford multiplication as in the introduction. A solution of $D\Psi=0$ is called a \textit{monogenic spinor}. See that the formula (\ref{k-Dirac operator}) differs from (\ref{Dirac operator locally}) only by replacing each coordinate vector field $\partial_{x_{\alpha i}}$ by $L_{\alpha i}$.

A real analytic $\Sp$-valued function $\Psi$ on $\aff$ can be written in a unique way as a converging sum $\sum_{i\ge0}\Psi_i$ where  $\wdeg(\Psi_i)=i$. As $D$ is a linear combination of vector fields $L_X$ with $X\in\lag_{-1}$, we see that $\wdeg(D\Psi_i)=i-1$. It follows that   $\Psi$ is a monogenic spinor iff  each $\Psi_i$ is a monogenic  spinor.  We will denote the set of all real analytic monogenic spinors over $\aff$ by $\ms$ and by $\ms_i$ the vector space all homogeneous monogenic spinors of the weighted degree $i$.  The map in (\ref{weighted jets of spinors at the origin}) then restricts to isomorphism of vector spaces 
\begin{eqnarray}\label{weighted monogenic spinors and tableau}
&\ms_i\ra\tab^{(i-1)}&\\
&\Psi\mapsto\xj^i_{x_0}\Psi.&
\end{eqnarray}  
In this paper we will be interested in the vector spaces $\ms_i$. Nevertheless, as we have just pointed out, this provides a good deal of information also about $\ms$. But we will not discuss this issue here (see \cite{MoIII} about the issue of the convergence of formal solutions). Let us recall the main result of this article which has been already stated in the introduction.

\begin{thm}
Let $i$ be a non-negative integer  and $\psi$ be a $\Sp$-valued function defined on the subset $M(n,k,\R):=\{x_{n+1,1}=\ldots=x_{n+1,k}=y_{12}\ldots=y_{k-1,k}=0\}$ of $\aff$ such that each component of $\psi$ is  a homogeneous polynomial of the degree $i$. Then there is a unique monogenic spinor $\Psi\in\ms_i$ such that $\Psi|_{M(n,k,\R)}=\psi$.
\end{thm}
Prove  of the uniqueness of $\Psi$:  Suppose that the claim is not true. Then  there is a non-constant monogenic spinor $\Psi$ which depends only on the variables $y_{rs},x_{n+1,i}$ where $i,r,s=1,\dots,k$. 
Then for any $1\le i<j\le k$ we have that $$0=\partial_{x_{1 i}}D_j\Psi=\partial_{x_{1 i}}\sum_{\beta=1}^{n+1} \varepsilon_\beta. L_{\beta j}\Psi=-\frac{1}{2}\varepsilon_1.\partial_{y_{ij}}\Psi.$$ We see that $\Psi$ does not depend also on any of the variable $y_{ij}$. So for any $i=1,\dots,k:0=D_i\Psi=\varepsilon_{n+1}.\partial_{x_{n+1,i}}\Psi$. We see that $\Psi$ is constant.  $\hfill\Box$ 

As we now know that the restriction map $\Psi\in\ms_i\mapsto\Psi|_{M(n,k,\R)}$ is injective, to finish the proof  it suffices to show that the dimension of $\ms_i$ is equal to the space of homogeneous spinors on $M(n,k,\R)$ of the weighted degree $i$. This is obviously true if $i=0$. From (\ref{dimension of tableaux}) follows that  this is true also if $i=1,2$. It remains to show the claim for $i=3,4,5,\dots$. This will occupy the rest of the paper.

\section{The Spencer complex for left invariant differential operators}
The Spencer complex (as  already introduced in \cite{Ta} and used in \cite{MoIII}, etc.) is a  complex of vector bundles  which is natural to filtered manifolds.  We will consider this complex over the origin. We will use a slightly different definition of the complex. More precisely we will not change the spaces in the complex but we will only change the definition of the co-differentials. Here we will take the co-differential with respect to the right invariant structure on $\aff$ rather than the natural left invariant structure. Then we get the Spencer complex associated to the $k$-Dirac operator which has the same form as the  Spencer complex used in the theory of exterior differential systems (see \cite{Br} and \cite{Sp}). In the second part of the section we define filtration on the tableau, on its prolongations and on all other spaces in the Spencer complex associated to the $k$-Dirac operator. As we mentioned and explained in the introduction, it is most natural to do this with respect to infinitesimal symmetries of the parabolic structure. Then the co-differential in the Spencer complex associated to the $k$-Dirac operator is nicely adapted to the filtration  and we will be able in the proof of Lemma \ref{lemma surjectivity} to repeat the proof of Proposition 2.5 from \cite{Br} which is the key point in the proof of Theorem \ref{main thm}.
  
Let us now recall the construction of the Spencer complex on  filtered manifolds. Let $i\in\N$. The exterior derivative of vector valued functions induces an injective linear map of vector spaces
\begin{eqnarray}\label{ex der of jets}
&\partial:\gr^{i+1}\Sp_\lambda\ra \lag_1\otimes\gr^{i}\Sp_\lambda\oplus\lag_2\otimes\gr^{i-1}\Sp_\lambda&\\
&\partial(\xj^{i+1}_{x_0}f)=(X\mapsto\xj^{i}_{x_0} (L_Xf),Y\mapsto\xj^{i-1}_{x_0}(L_Yf)).&\nonumber
\end{eqnarray}
where $X\in\lag_{-1},Y\in\lag_{-2}$ and $f$ is the germ of $\Sp_\lambda$-valued function at $x_0$ which satisfies $\xj^i_{x_0}f=0$. We have used here the notation from (\ref{notation for fibres of weighted jets over the origin for homogeneous bundles}), isomorphisms $\lag_i\cong\lag_{-i}^\ast,i=1,2$ and we agree that $\gr^{i}\Sp_\lambda=0$ if $i<0$. Notice that the map (\ref{ex der of jets}) is well defined at any point $x\in\aff$ and that even the first component of $\partial$ is  injective   (as $H$ is bracket generating). More generally, the exterior derivative of $\Sp_\lambda$-valued $r$-forms gives linear map
\begin{eqnarray}
\partial:\Lambda^r\gr^{i+1-r}\Sp_\lambda\ra
\Lambda^{r+1}\gr^{i-r}\Sp_\lambda
\end{eqnarray}
where we put
\begin{equation}
\Lambda^\bullet\gr^\ast\Sp_\lambda:=\bigoplus_{\ell=0}^{\lfloor\frac{\bullet}{2}\rfloor}\Lambda^{\bullet-\ell}\lag_1\wedge\Lambda^{\ell}\lag_2\otimes\gr^{\ast-\ell}\Sp_\lambda.
\end{equation}
We obtain a complex $(\Lambda^\ast\gr^{i+1-\ast}\Sp_\lambda,\partial)$ with trivial cohomology groups (this follows from the fact that the de Rham complex is locally exact). Well known properties of the exterior derivative imply that 
\begin{eqnarray}\label{Leibniz rule}
&\partial(\omega\wedge\omega')=\partial\omega\wedge\omega'+ (-1)^\bullet\omega\wedge\partial\omega',&\\
&\partial(e^i\otimes\varepsilon_\alpha)=0,\ \partial(e^r\wedge e^s)=\sum_{\alpha=1}^{n+1}(e^r\otimes\varepsilon_\alpha)\wedge (e^s\otimes\varepsilon_\alpha).&\label{kostant differential}
\end{eqnarray}
where $\omega\in\Lambda^\bullet\gr^\ast\Sp_\lambda,\omega'\in\Lambda^{\bullet'}\gr^{\ast'}\Sp_\lambda$. The line (\ref{kostant differential}) follows from  (\ref{ex of left inv contact form}). 
 However,  as we already mentioned in  Remark \ref{remark differentiating monogenic spinors}, it is not true that if we restrict  $\partial$  in (\ref{ex der of jets})  to the subspace $\tab^{(i)}$ that we get a map $\tab^{(i)}\ra\lag_1\otimes\tab^{(i-1)}\oplus\lag_2\otimes\tab^{(i-2)}$. The problem is that the left invariant vector fields are not infinitesimal symmetries of the structure and so they do not commute with the $k$-Dirac operator $D$. To fix this problem we need to consider exterior derivative  with respect to right invariant objects on $\aff$ rather than left invariant objects. This leads us to consider the map:
\begin{eqnarray}\label{right invariant differential}
&\delta:\gr^{i+1}\Sp_\lambda\ra \lag_1\otimes\gr^{i}\Sp_\lambda\oplus\lag_2\otimes\gr^{i-1}\Sp_\lambda&\\
&\delta(\xj^{i+1}_{x_0}f):=(X\mapsto\xj^{i}_{x_0} (R_Xf),Y\mapsto\xj^{i-1}_{x_0}(R_Yf)).&
\end{eqnarray}
where  $X\in\lag_{-1},Y\in\lag_{-2},\xj^{i+1}_{x_0}f$ are as in (\ref{ex der of jets}) and more generally:
\begin{eqnarray}\label{sc}
\delta:\Lambda^{r}\gr^{i+1-r}\Sp_\lambda\ra\Lambda^{r+1}\gr^{i-r}\Sp_\lambda
\end{eqnarray}
such that (\ref{Leibniz rule}) is still true with $\partial$ being replaced by $\delta$, but from (\ref{differential of contact forms}) follows that (\ref{kostant differential}) has to be replaced by
\begin{eqnarray}\label{differential of forms}
\delta(e^i\otimes\varepsilon_\alpha)=0,\ \delta(e^r\wedge e^s)=-\sum_{\alpha=1}^{n+1}(e^r\otimes\varepsilon_\alpha)\wedge (e^s\otimes\varepsilon_\alpha).
\end{eqnarray}
 This uniquely pins down $\delta$. Notice that we have not changed the spaces in the complexes but we have only replaced the map $\partial$ by the map $\delta$. By the same reason as above, the cohomology of the complex $(\Lambda^\ast\gr^{i+1-\ast}\Sp_\lambda,\delta)$ is  trivial. Now, we can restrict the map (\ref{right invariant differential}) to the subspace $\tab^{(i)}$ and we get an injective map
\begin{equation}
\tab^{(i)}\ra\lag_1\otimes\tab^{(i-1)}\oplus\lag_2\otimes\tab^{(i-2)}
\end{equation}
which we still denote by $\delta$. More generally, we can restrict (\ref{sc}) also to 
\begin{equation}
\tab^{\bullet,\ast}:=\bigoplus_{\ell=0}^{\lfloor\frac{\bullet}{2}\rfloor}\Lambda^{\bullet-\ell}\lag_1\wedge\Lambda^{\ell}\lag_2\otimes\tab^{(\ast-\ell)}.
\end{equation}
It follows that  the  complex   $(\Lambda^\ast\gr^{i+1-\ast}\Sp_\lambda,\delta)$ contains the subcomplex $(\tab^{\ast,i-\ast},\delta)$ which is the \textit{Spencer complex associated to the $k$-Dirac operator.}

\subsection{Filtration of the tableau}
The map (\ref{differentiating with right invariant vector fields}) is a special case of the (linear map induced by the) Lie derivative $\cL_{R_X}:\Lambda^\bullet\gr^{\ast}\Sp_\lambda\ra\Lambda^\bullet\gr^{\ast-1}\Sp_\lambda$. This  map  satisfies $\cL_{R_X}(\omega\wedge\omega')=(\cL_{R_X}\omega)\wedge\omega'+\omega\wedge(\cL_{R_X}\omega')$ where $\omega,\omega'$ are as in (\ref{Leibniz rule}). 
The Cartan formula $\cL_{R_X}=i_{X}\delta+\delta i_{X}$ is still valid where $i_{X}$ is the insertion of $X\in\lag_-$ into the first entry. Notice that from (\ref{differential of forms}) follows
\begin{eqnarray}\label{Lie derivative of forms}
&\cL_{R_{\alpha i}} (e^j\otimes\varepsilon_\beta)=0,\ \cL_{R_{\alpha i}}(e^r\wedge e^s)=\delta_i^s e^r\otimes\varepsilon_\alpha-\delta_i^r e^s\otimes\varepsilon_\alpha,&\\
&\cL_{\partial_{y_{rs}}} (e^j\otimes\varepsilon_\beta)=\cL_{\partial_{y_{rs}}}(e^u\wedge e^v)=0.&\label{Lie derivative by transversal symmetry}
\end{eqnarray}
The Lie derivative restricts to a map
\begin{equation}
\cL_{R_X}:\tab^{\bullet,\ast}\ra\tab^{\bullet,\ast-1}.
\end{equation}

Now we will introduce a filtration on the spaces $\tab^{\bullet,\ast}$. Let us fix a basis $\{X_1,\dots,X_{k(n+1)}\}$ of $\lag_{-1}$. We will for simplicity write $R_p$ instead of $R_{X_p},p=1,\dots,k(n+1)$. Then we put for each $r,i$ and $j=0,\dots,k(n+1)$:
\begin{equation}
\tab^{r,i}_j:=\{s\in\tab^{i,r}|\cL_{R_{1}}s=\dots=\cL_{R_{j}}s=0\}.
\end{equation}
 We obtain a filtration 
\begin{equation}\label{filtration}
\{0\}=\tab^{r,i}_{k(n+1)}\subset\tab^{r,i}_{k(n+1)-1}\subset\dots\subset\tab^{r,i}_1\subset\tab^{r,i}_0=\tab^{r,i}
\end{equation}
 of $\tab^{r,i}$.  We will for simplicity write $\tab^{(i)}_j$ instead of $\tab^{0,i}_j$.
 
 We have now given all necessary definitions. Now we can proceed with the proof of Theorem \ref{main thm}. Let us briefly go through the next steps. We will start with  Lemmas \ref{lemma help} and  \ref{lemma preferred basis} which are needed in order to make the machinery of the  Cartan-K\"ahler theorem  running also for the $k$-Dirac operator. Notice that in contrast to the case of the classical Cartan-K\"ahler theorem, it is not at all clear that Lemma \ref{lemma preferred basis} is true. The problem here is that the right invariant vector fields do not commute but rather:
\begin{equation}\label{commutator of right invariant fields}
[\cL_{R_{\alpha u}},\cL_{R_{\beta v}}]=-\delta_{\alpha\beta}\cL_{\partial_{y_{uv}}}.
\end{equation}
However, we will show that if we choose the basis 
\begin{eqnarray}\label{preferred basis}
&\{e_1\otimes\varepsilon_1,\dots,e_1\otimes\varepsilon_n,e_2\otimes\varepsilon_1,\dots,e_2\otimes\varepsilon_n,\dots,e_k\otimes\varepsilon_1,&\\
&\dots,e_k\otimes\varepsilon_n,e_1\otimes\varepsilon_{n+1},\dots,e_k\otimes\varepsilon_{n+1}\},&\nonumber
\end{eqnarray}
then everything works just as for the classical Cartan-K\"ahler theorem.  Then it follows that  for each $i,j$ there is a sequence
\begin{equation}\label{ses}
0\ra\tab^{(i+1)}_{j}\ra\tab^{(i+1)}_{j-1}\xrightarrow{R_{j}}\tab^{(i)}_{j-1}\ra0.
\end{equation} 
Clearly, the sequence is a complex which is exact in the middle by the definition of $\tab^{(i+1)}_{j}$. As $\tab^{(i+1)}_{k(n+1)}=0$, we obtain  an upper bound on the dimension of $\tab^{(i+1)}=\tab^{(i+1)}_0$:
\begin{equation}\label{upper bound}
 \dim\tab^{(i+1)}\le \sum_{j=0}^{k(n+1)}\dim\tab^{(i)}_j.
\end{equation} 
 Notice that there is  equality  in (\ref{upper bound})   iff  $R_j$ is surjective for each $j=1,2,\dots,(n+1)k$.  We will show in Lemma \ref{lemma surjectivity} that the latter condition holds. In other words, we will prove that (\ref{ses}) is a short exact sequence  (in the language of the classical Cartan-K\"ahler theorem this property is the definition of the involutivity of the system). 
The fact that (\ref{ses}) is a short exact sequence for each $i,j$ is all we need to  finish the proof of Theorem \ref{main thm} (see the end of the section).
After this short summary, we can proceed by verifying all steps.  
\begin{lemma}\label{lemma help}
Let $i\ge0$ be an  integer and  $\{\tab^{(i+1)}_j:j=0,\dots,k(n+1)\}$ be the filtration of $\tab^{(i+1)}$ with respect to the basis (\ref{preferred basis}).
  If $f\in\tab^{(i+1)}_{nr}$ where $r=1,\dots,k-1$, then $\partial_{y_{st}}f=0\in\tab^{(i-1)}$ whenever $s\le r$.
\end{lemma}

\begin{remark}\label{remark monogenic spinors in less variables}
Before going through the proof, let us make few observations. We will view the element $f\in\tab^{(i+1)}_{nr}$ from the statement of Lemma \ref{lemma help} as a homogeneous monogenic spinor of the weighted degree $i+2$ (here we use the isomorphism from  (\ref{weighted monogenic spinors and tableau})). If $f$ satisfies the hypothesis and the conclusion of the lemma, then $f$ depends only on the variables $x_{\alpha j},y_{cd}$ where $j>r,d>c>r$ and $\alpha=1,\dots,n+1$. This easily follows from the following computation. Suppose that $s\le r,\alpha=1,\dots,n$. Then we have
\begin{eqnarray}
&&0=R_{\alpha s}f=\partial_{x_{\alpha s}}f+\frac{1}{2}\sum_{t=1}^kx_{\alpha t}\partial_{y_{st}}f=\partial_{x_{\alpha s}}f.\\
&&0=D_sf=\sum_{\beta=1}^{n+1}\varepsilon_\beta.L_{\beta s}f=\sum_{\beta=1}^{n+1}\varepsilon_{\beta}.\partial_{x_{\beta s}}f=\varepsilon_{n+1}.\partial_{x_{n+1,s}}f.
\end{eqnarray}
In particular, if $f\in\tab^{(i+1)}_{n(k-1)}$, then $f$ depends only on the variables $x_{\alpha k},\alpha=1,\dots,n+1$. Then $D_if=0,i=1,\dots,k-1$ and $D_kf$ is the usual Dirac operator (in one variable) which  is an involutive system.
It also follows from (\ref{Lie derivative of forms}) that 
\begin{equation}\label{filtration in less variables}
\tab^{r,i}_{j}\subset\bigoplus_{\ell=0}^{\lfloor\frac{r}{2}\rfloor}\Lambda^{r-\ell}\lag_1\wedge\Lambda^{\ell}\lag_2\otimes\tab^{(i-\ell)}_{nJ}
\end{equation}
where we write $j=nJ+\rho,\rho=1,\dots,n-1$.
\end{remark} 

Proof of Lemma \ref{lemma help}:   Let us suppose that the claim is not true. Then we may choose $f$ (which we view as an element of $\ms_{i+2}$ as in Remark \ref{remark monogenic spinors in less variables}) which satisfy the hypothesis but $\partial_{y_{st}}f\ne0$ for some $s\le r,t>s$. We will show that this leads to a contradiction. We may also assume that $i,r$ are minimal,  i.e. if  $g\in\tab^{(t+1)}_{nj}$ satisfies $\partial_{y_{st}}g\ne0$  for some $s\le j$, then  $t>i$ or $t=i,j\ge r$. In particular, as $f\in A^{(i+1)}_{nr}\subset A^{(i+1)}_{ns}$, we have that  $s=r$ by the choice of $f$. 
As we have showed in Remark \ref{remark monogenic spinors in less variables},  $f$ may depend only on the variables $x_{\alpha j},y_{cd}$ where $j\ge r,d>c\ge r,\alpha=1,\dots,n+1$. 

We have that  $f_{uv}:=\partial_{y_{uv}}f\in\tab^{(i-1)}_{nr}$. It follows from the choice of $f$ and by Remark \ref{remark monogenic spinors in less variables} that $f_{uv}$ depends only the variables $x_{\alpha j},y_{cd}$ where $j> r,d>c> r,\alpha=1,\dots,n+1$.  This implies that 
\begin{equation}\label{x}
f=\sum_{u>r}y_{ru}f_{ru}+\dots
\end{equation}
 where  $\dots$ is  a function which does not depend on any of the variable $y_{cd},c\le r$. By the assumption above, $f_{rt}\ne0$.  From (\ref{Lie bracket}) follows that $R_{n+1,j}f\in\tab^{(i)}_{nr};j=1,\dots,k$ and hence as before,  $R_{n+1,j}f$ does not depend on any of the variable $y_{cd},c\le r$. Using (\ref{x}), we have 
\begin{equation}
R_{n+1,j}f=\sum_{u>r}y_{ru}\partial_{x_{n+1,j}}f_{ru}+\dots
\end{equation}
 where  $\dots$ is  a function  which does not depend on any of the variable $y_{ru},u=1,\dots,k$. It follows that $f_{ru}$ does not depend on $x_{n+1,j};j=1,\dots,k$. We have that 

\begin{eqnarray}
&&0=R_{\alpha r}f=\sum_{u>r}\frac{1}{2}x_{\alpha u}\partial_{y_{ru}}f+\partial_{x_{\alpha r}}f=\sum_{u>r}\frac{1}{2}x_{\alpha u}f_{ru}+\partial_{x_{\alpha r}}f;\alpha<n\label{z}\\
&&0=\sum_{\alpha=1}^{n+1}\varepsilon_\alpha.L_{\alpha r}f=\sum_{\alpha\le n}\varepsilon_\alpha.(R_{\alpha r}-\sum_{u>r}x_{\alpha u}\partial_{y_{ru}})f+\varepsilon_{n+1}.L_{n+1,r}f\nonumber\\
&&=-\sum_{\alpha\le n,u>r}\varepsilon_\alpha.(x_{\alpha u}f_{ru})+\varepsilon_{n+1}.L_{n+1,r}f.\label{y}
\end{eqnarray}
Now the first term in (\ref{y}) depends only the variables $x_{\alpha j},y_{cd}$ where $j> r,d>c>r,\alpha=1,\dots,n$. In particular, it does not depend on $x_{n+1,j};j=1,\dots,k$. This  implies that the same is true for
$$L_{n+1,r}f=\varepsilon_{n+1}.(\partial_{x_{n+1,r}}-\frac{1}{2}\sum_{u>r}x_{n+1,u}\partial_{y_{ru}})f=\varepsilon_{n+1}.(\partial_{x_{n+1,r}}f-\frac{1}{2}\sum_{u>r}x_{n+1,u}f_{ru}).$$
Combining this together with (\ref{z}), we obtain that
$$f=\sum_{u>r}(y_{ru}f_{ru}-\frac{1}{2}\sum_{\alpha=1}^nx_{\alpha r}x_{\alpha u}f_{ru}+\frac{1}{2}x_{n+1,r}x_{n+1,u}f_{ru})+\dots$$
where $\dots$ is a function  which depends only on $x_{\alpha c},y_{uv}$ where $c,u,v>r$. Then
\begin{eqnarray*}
0&=&\sum_{\alpha=1}^{n+1}\varepsilon_\alpha. L_{\alpha t}f=\sum_{\alpha\le n}\varepsilon_\alpha.L_{\alpha t}f+\varepsilon_{n+1}.L_{n+1,t}f\\
&=&\sum_{\alpha\le n}\varepsilon_\alpha.(-\frac{1}{2}x_{\alpha r}f_{rt}+\frac{1}{2}x_{\alpha r}f_{rt})+\frac{1}{2}\varepsilon_{n+1}.(x_{n+1,r}f_{rt}+x_{n+1,r}f_{rt})+\dots
\end{eqnarray*}
where $\dots$ represents a function which depends on $x_{\alpha c},y_{uv}$ where $c,u,v>r,\alpha=1,\dots,n+1$. We see that the equality holds iff  $f_{rt}=0$. The claim is proved.  $\hfill\Box$\\

We will formulate Lemma \ref{lemma surjectivity} not  only for the basis (\ref{preferred basis}) but for a particular set of bases of $\lag_{-1}$ although this is not necessary in order to show that (\ref{ses}) is a well defined sequence. We will need this more general statement in  the proof of Lemma \ref{lemma surjectivity}.

\begin{lemma}\label{lemma preferred basis}
Let $\sigma$ be a permutation of $\{1,2,\dots,n\}$. Let 
$\{\tab^{r,i}_j[\sigma]:j=0,\dots,(n+1)k\}$ be the filtration of $\tab^{r,i}$ with respect to the basis
\begin{eqnarray}\label{permuted basis}
&\{e_1\otimes\varepsilon_{\sigma(1)},\dots,e_1\otimes\varepsilon_{\sigma(n)},e_2\otimes\varepsilon_{\sigma(1)},\dots,e_2\otimes\varepsilon_{\sigma(n)},\dots,e_k\otimes\varepsilon_{\sigma(1)},&\\
&\dots,e_k\otimes\varepsilon_{\sigma(n)},e_1\otimes\varepsilon_{n+1},\dots,e_k\otimes\varepsilon_{n+1}\}&\nonumber
\end{eqnarray}
 of $\lag_{-1}$. Then for each $i,r,j:\cL_{R_j}f\in\tab^{r,i+1}_{j-1}[\sigma]$ whenever $f\in\tab^{r,i}_{j-1}[\sigma]$. 
\end{lemma}
Proof:   The claim follows from (\ref{commutator of right invariant fields}),  the Leibniz property of  Lie derivative, the formulas (\ref{Lie derivative by transversal symmetry}),  (\ref{filtration in less variables}) and Lemma \ref{lemma help}. $\hfill\Box$\\

We will use the following notation. Let us for simplicity denote the  basis (\ref{preferred basis}) by $\{X_1,\dots,X_{k(n+1)}\}$  and let $\{\varpi_1,\dots,\varpi_{k(n+1)}\}$ be the dual  basis of $\lag_1$  so that $\varpi_p(X_q)=\delta_{pq}$ for each $p,q$. 
We will  put $e^{rs}:=e^r\wedge e^s$ so we obtain a basis  $\{e^{12},\dots,e^{k-1,k}\}$ which is dual to the basis $\{e_1\wedge e_2,\dots,e_{k-1}\wedge e_k\}$ of $\lag_{-2}$.  Then we can write down the map (\ref{right invariant differential}) as

\begin{equation}\label{differential with framing}
\delta(f)=\sum_{l=1}^{k(n+1)}\varpi_l\otimes R_lf+\sum_{1\le r<s\le k}e^{rs}\otimes\partial_{y_{rs}}f
\end{equation}
where $f\in\tab^{(i)}$.

By Remark \ref{remark monogenic spinors in less variables} we have that  $\dim(\tab^{(i)}_j)=0$ for $j\ge kn$. In Section 4.3 in \cite{S} can be found  that 
\begin{equation}\label{cartan characters}
\dim(\tab^{(0)}_j)=
\bigg\{
\begin{matrix}
(nk-j)\dim\Sp; \ j=0,\dots,nk,\\
0;\ \ \ \ j\ge nk.\\
\end{matrix}
\end{equation}
It can be easily seen that the formula holds also if we replace $\tab^{(0)}_j$ by $\tab^{(0)}_j[\sigma]:=\tab_j^{0,0}[\sigma]$ where $\sigma$ is any  permutation from Lemma \ref{lemma preferred basis}.  As $\dim(\tab^{(1)})={nk+1\choose2}\dim\Sp$ (see formula (\ref{dimension of tableaux})), we have that
\begin{equation}
\dim\tab^{(1)}=\sum_{j=0}^{k(n+1)}\dim(\tab^{(0)}_j[\sigma])
\end{equation}
and so we have that for each $j=1,2,\dots,k(n+1)$ the map 
\begin{equation}\label{surjectionI}
R_j:\tab^{(i+1)}_{j-1}[\sigma]\ra\tab^{(i)}_{j-1}[\sigma]
\end{equation}
is surjective if $i=0$. 

\begin{lemma}\label{lemma surjectivity}
Let $\sigma$ be a permutation of $\{1,\dots,n\},i\ge0$ and $j=1,\dots,k(n+1)$. Then the map $R_j$ in (\ref{surjectionI}) is surjective.
 \end{lemma}
Proof: As we have just seen above, the claim is true if $i=0$. Let $r\ge0,p=1,\dots,k(n+1)$ and assume that (\ref{surjectionI}) is surjective whenever $\sigma$ is a permutation of $\{1,\dots,n\}$ and $i<r,j=1,\dots,k(n+1)$ or $i=r,j=p+1,\dots,k(n+1)$. We will show that then (\ref{surjectionI}) is surjective also when $\sigma$ is the identity permutation and $i=r,p=j$. The proof for arbitrary permutation is similar. Let us consider the following commutative diagram 

\begin{equation}
\xymatrix{\tab^{(r+1)}_{p-1}\ar[r]^\delta\ar[d]^{R_p}&\tab^{1,r}_{p-1}\ar[r]^\delta\ar[d]^{\cL_{R_p}}&
\tab^{2,r-1}_{p-1}
\ar[d]^{\cL_{R_p}}\ar[r]^\delta&\dots\\
\tab^{(r)}_{p-1}\ar[r]^\delta&\tab^{1,r-1}_{p-1}\ar[r]^\delta&\tab^{2,r-2}_{p-1}\ar[r]^\delta&\dots}
\end{equation}
Let $Q\in \tab^{(r)}_{p-1}$. Suppose that there is  $T\in\tab^{1,r}_{p-1}$ such that the following is true:
\begin{enumerate}
\item $R_pT=\delta Q$.
\item $\delta T=0$.
\item $i_{R_1}T=\dots=i_{R_{p-1}}T=0$.
\end{enumerate}
 Then by the exactness of the complex (\ref{sc}) there is $P\in\gr^{r+2}\Sp_\lambda$ such that $\delta P=T$.  As  $R_sP=i_{R_s}\delta P=i_{R_s}T;s=1,\dots,k(n+1)$ we have that $R_1 P=\dots=R_{p-1}P=0$. Since $R_sP\in\tab^{(r)}$ for each $s=1,\dots,k(n+1)$ it follows that $P\in\tab^{(r+1)}_{p-1}$. Then $\delta R_p(P)=\cL_{R_p}\delta(P)=\cL_{R_p}(T)=\delta Q$. By the injectivity of $\delta$ at the beginning, it follows that $Q= R_p(P)$. Hence,  it is enough to find $T$.

If $p> nk$, then by  Remark \ref{remark monogenic spinors in less variables} we have  $\tab^{(r+1)}_{p-1}=\tab^{(r)}_{p-1}=0$  and so there is nothing to prove. If $nk\ge p> k(n-1)$, then by the Remark \ref{remark monogenic spinors in less variables} we are dealing with the Dirac operator, this is an involutive system and so the claim follows also in this case. So we may actually assume that $p\le n(k-1)$. Write $p=nJ+\rho$ where $0<\rho\le n$. Then $J=0,1,\dots,k-2$ and $X_p=e_{J+1}\otimes\varepsilon_\rho$. By Lemma \ref{lemma help} we have that $\delta Q=\sum_{l\ge p}\varpi_{l}\otimes q_{l}+\sum_{J< c<d}e^{cd}\otimes q_{cd}$ where $q_{l}=R_{X_{l}}Q,q_{rs}=\partial_{y_{cd}}Q$. We will first show that there is  $U\in\tab^{1,r}_{p-1}$ 
such that $\cL_{R_p}U=\delta Q$. If  $q_{cd}=0$ for all $c,d$, then $q_l\in\tab^{(r-1)}_{p-1}$ for each $l\ge p$ and so by the induction hypothesis  there are $u_l\in\tab^{(r)}_{p-1},l\ge p$ such that $R_pu_l=q_l$ and so we may put $U=\sum_{l\ge p}\varpi\otimes u_l$.
In general, each $q_{cd}\in\tab^{(r-2)}_{p-1}$  and so by the induction hypothesis there are $u_{cd}\in\tab^{(r-1)}_{p-1}$ such that $R_{p}u_{cd}=q_{cd}$.  Let us fix $t=1,\dots\rho-1$ and let $\sigma_t$ be the transposition $(t,\rho-1)$. We have that $\tab^{(r-1)}_{p-1}=\tab^{(r-1)}_{p-1}[\sigma_t]\subset\tab^{(r-1)}_{p-2}[\sigma_t]$. So by the induction hypothesis:  for each $d>J+1,t=1,\dots,p-1$ there is $u'_{dt}\in \tab^{(r)}_{p-2}[\sigma_t]$ such that $R_{nJ+t}u'_{dt}=u_{J+1,d}$. Notice that $X_{nJ+t}=e_{J+1}\otimes\varepsilon_t$ and so from (\ref{Lie derivative of forms}) follows that  $\cL_{R_{nJ+t}}e^{cd}=-\delta^{c}_{J+1}\varpi_{(d-1)n+t}$ whenever $J<c<d$.  Consider  
$$U':=\sum_{J+1<d,t=1,\dots,\rho-1}\varpi_{n(d-1)+t}\otimes u'_{dt}+\sum_{J<c<d}e^{cd}\otimes u_{cd}.$$
 Then we clearly have that $U'\in\tab^{1,r}_{nJ}$. For each $t=1,\dots,\rho-1$ we have  $$\cL_{R_{nJ+t}}U'=\sum_{J+1<d}\varpi_{n(d-1)+t}\otimes u_{J+1,d}-\varpi_{(d-1)n+t}\otimes u_{J+1,d}=0$$
  and so we see actually see that $U'\in\tab^{(r)}_{p-1}$. Then we have already proved that there is  $U''\in\tab^{r,1}_{p-1}$ such that $\cL_{R_p}U''=\delta Q-\cL_{R_p}U'$ and so we can put $U:=U'+U''$. 



 We have that $V:=\delta U\in\tab^{2,r-1}_{p-1}$ and  for each $j<p:i_{R_j}\delta U=\cL_{R_j} U=0$. Since $\cL_{R_p}\delta U=\delta \cL_{R_p}U=\delta^2Q=0$
we find that  $V\in \tab^{2,r-1}_p$. As $U\in\tab^{1,r}_{nJ}$ it follows by (\ref{filtration in less variables}) that 

\begin{eqnarray*}
&&V:=\varpi_p\wedge(\sum_{l>p}\varpi_l\otimes v'_{l}+\sum_{J<c<d}e^{cd}\otimes v'_{cd})
+\sum_{b>a>p}\varpi_a\wedge\varpi_b\otimes v''_{ab}\\
&&+\sum_{a>p,J<c<d}\varpi_a\wedge e^{cd}\otimes v''_{acd}+\sum_{J<c<d,J<u<v}e^{cd}\wedge e^{uv}\otimes v''_{cduv}.
\end{eqnarray*}
 Now $\delta i_{R_p}\delta U=\delta \cL_{R_p}U=0$ and so  $V':=i_{R_p}\delta U$ and $V'':=V-\varpi_p\wedge V'$ are two closed forms. 
   It immediately follows from  (\ref{differential of forms}), (\ref{differential with framing})  that  $v'_{p+1}\in\tab^{(r-1)}_{p},v''_{p+1,b}\in\tab^{(r-1)}_p,b\ge p+2$. So we may find $t_p\in\tab^{(r)}_{p}$ such that $\cL_{R_{p+1}}t_p=v'_p$ and so $V_1':=V'-\delta t_p=\sum_{l\ge p+1}\varpi_l\otimes v'_{1l}+\sum_{J<c<d}e^{cd}\otimes v'_{1cd}$. Repeating the  argument for $V_1'$ we eventually construct\footnote{Here we also use the following observation: if the 1-form $V'$ was closed such that $v'_{p+1}=\dots=v'_{k(n+1)}=0$, then $V'=0$.}   $T_1\in\tab^{(r-1)}_p$ such that $\delta T_1=V'$.  Then $\delta(-\varpi_p\otimes T_1)=\varpi_p\wedge V'$. Similarly, there are  $t_{p+1,b}\in\tab^{(r)}_p$ such that  $\cL_{R_{p+1}}t_{p+1,b}=v''_{p+1,b},b\ge p+2$. It follows that 
   $$V''_1:=V''-\delta(\sum_{b\ge p+2}\omega_{b}\otimes t_{p+1,b})=\sum_{a>b> p+1}\varpi_a\wedge\varpi_b\otimes v''_{1ab}+\dots$$ We may again repeat the argument for $V''_1$ and we eventually construct  $T_2\in\tab^{(1,r)}_{p}$ such that $\delta T_2=V''$. So  $T:=U-T_1-T_2$ is a form we are looking for. $\hfill\Box$\\
   
Now we can complete the proof of the main theorem.\\

Proof of Theorem \ref{main thm}: As we have already proved injectivity of the restriction map $\Psi\in\ms_{i+1}\mapsto\Psi|_{M(n,k,\R)}$, we know that $\dim\ms_{i+1}\le{nk+i\choose i+1}\dim\Sp$. Hence, it is enough to show that there is equality. 
We will prove by induction on $p,r$ that 
\begin{equation}\label{xy}
\dim\tab^{(p)}_r={nk-r+p\choose p+1}\dim\Sp.
\end{equation}
 Check that the right hand side is equal to the dimension of the space of homogeneous polynomials of the degree $p+1$ in $nk-r$ variables multiplied by $\dim\Sp$ (where we agree that the binomial coefficient is zero if $nk-r+p<p+1$).
 By (\ref{cartan characters}) the claim is true for $p=0,r=0,1,\dots,(n+1)k$. Let us fix non-negative integers $i,j$.  We suppose that (\ref{xy}) is true  if  $p<i$ or $p=i,k(n+1)\ge r>j$.
  We want to prove the claim also  $i,j$. As we already know that  (\ref{ses}) is short exact for each $i,j$, we have  that
\begin{eqnarray*}
\dim\tab^{(i)}_{j}&=&\dim\tab^{(i)}_{j+1}+\dim\tab^{(i-1)}_{j}\\
&=&\bigg({nk-j-1+i\choose i+1}+{nk-j+i-1\choose i}\bigg)\dim\Sp\\
&=&{nk-j+i\choose i+1}\dim\Sp.
\end{eqnarray*}
The equality we wanted to prove is the particular case $j=0$ as $\dim\ms_{i+1}=\dim\tab^{(i)}$. $\hfill\Box$

\def\bibname{Bibliography}


\begin{thebibliography}{99}
\addcontentsline{toc}{chapter}{\bibname}


\bibitem{Br} {\sc Bryant} R. L. , S. S. {\sc Chern}, R. B. {\sc Gardner}, H. L. {\sc Goldschmidt} and P. A. {\sc Griffiths}. \emph{Exterior differential systems, Mathematical Sciences Research Institute Publications}, vol. 18, Springer-Verlag, New York, 1991. 





\bibitem{CS}
  {\sc \v Cap,} Andreas, Jan {\sc Slov\'ak}.
  \emph{Parabolic Geometries I, Background and General Theory}.
  American Mathematical Society, Providence, 2009.
ISBN 978-0-8218-2681-2.


\bibitem{CSSS}
  {\sc Colombo,} Fabrizio, Irene {\sc Sabadini}, Franciscus {\sc Sommen}, Daniele C. {\sc Struppa}.
  \emph{ Analysis of Dirac Systems and Computational Algebra}.
Birkhauser, Boston, 2004. ISBN 0-8176-4255-2.



\bibitem{GW} \sc{Goodman} Roe, Nolan R. \sc{ Wallach}: \emph{Representations and Invariants of the Classical Groups}. Springer, New York, 2009. ISBN 978-0-387-79581-6

\bibitem{K} \sc{Kruglikov} Boris.
\emph{Symmetries of  filtered structures via  filtered Lie equations}. Journal of Geometry and Physics, vol. 85, 2014, p. 164-170.


\bibitem{MoI}
{\sc Morimoto}, Tohru. \emph{Th\'eor\`eme de Cartan-K\"ahler dans une classe de fonctions formelles Gevrey.}
C. R. Acad. Sci. Paris. 311, s\`erie A. 1990, p. 443-436.

\bibitem{MoII}
{\sc Morimoto}, Tohru. \emph{Th\'eor\`eme d'existence de solutions analytiques pour des syst\`emes d'\'equations aux d\'eriv\'ees partielles non-lin\'eaires avec singularit\'es}. C.R. Acad. Sci. Paris. 321, s\'erie 1.
1995. p. 1491-1496.


\bibitem{MoIII}
{\sc Morimoto}, Tohru. \emph{Lie algebras, geometric structures and differential equations on filtered manifolds.}
In Lie Groups Geometric Structures and Differential Equations - One Hundred Years
after Sophus Lie, Adv. Stud. Pure Math., Math. Soc. of Japan, Tokyo. 2002, p. 205-252.


\bibitem{MoIV} 
{\sc Morimoto}, Tohru: \emph{Differential Equations Associated to a Representation of a Lie algebra from
the Viewpoint of Nilpotent Analysis}. RIMS Kokyuroku 1502, Kyoto University. 2006/07, p.
238-250.

\bibitem{MoV}{\sc Morimoto}, Tohru.
\emph{Generalized Spencer Cohomology Groups and
Quasi-Regular Bases.} Tokyo J. Math. vol. 14, no. 1, 1991, p. 165-179.

\bibitem{Ne} {\sc Neusser}, Katharina. \emph{ Prolongation on regular infinitesimal flag manifolds.} International Journal of Mathematics.
vol. 23, no. 4, 2012. p. 1-41.


\bibitem{S}{\sc Sala\v c}, Tom\'a\v s. \emph{k-Dirac operator and the Cartan-K\"ahler theorem}.
Archivum mathematicum. vol. 49, no. 5, 2013, p. 333-346.


\bibitem{Sa}{\sc Sala\v c}, Tom\'a\v s. \emph{k-Dirac operator and parabolic geometry}.
J.  comp. anal. and oper. theo. vol. 8, no. 2, 2014, p. 383-408.


\bibitem{SSSL}
  {\sc Sabadini,} Irene, Franciscus {\sc Sommen}, Daniell C. {\sc Struppa}, Peter van {\sc Lancker}.
  \emph{Complexes of Dirac operators in Clifford algebras}.
Mathematische Zeitschrift. vol. 239, no. 2, 2002, p. 293-320. 


\bibitem{SS} {\sc Slov\'{a}k,} Jan, Vladim\'ir {\sc Sou\v cek}.   \emph{Invariant Operators of the First Order on Manifolds with a Given Parabolic Structure}. Societ\'e Math\'ematique de France.  vol. 4, 2001, p. 251-276.

\bibitem{Sp}{\sc Spencer}, D.C. \emph{Overdetermined systems of linear partial differential equations}. Bull. Amer. Math. Soc., vol. 75, no. 2, 1969, p. 179-239.

\bibitem{Ta}{\sc Tanaka}, Noboru.
\emph{On the equivalence problems associated with simple graded Lie algebras}. Hokkaido Math.
J., vol. 8, no. 1,  1979, p. 23-84.
\end{thebibliography}
\end{document}